\documentclass[12pt]{article}
\usepackage{amsmath,amssymb,amsthm,eucal}
\usepackage[
  pdftitle={On the representation of bent functions by bent rectangles},
  pdfauthor={Sergey Agievich},
  bookmarks=true,pdfpagemode=None,pdfstartview=FitH]{hyperref}

\newcommand{\q}{\phantom{-}}
\newcommand{\func}[1]{{\mathcal #1}}
\newcommand{\F}{{\mathbb F}}

\newcommand{\tinybox}{\vrule\vbox to.333em{\hrule width.25em\vfil\hrule}\vrule}

\newcommand{\up}[2]{\mathop{\mathstrut #1}\limits^{#2}\mathstrut}
\newcommand{\circup}[1]{\up{#1}{\circ}}
\newcommand{\hatup}[1]{\up{#1}{\mbox{\tiny$\wedge$}}}
\newcommand{\boxup}[1]{\up{#1}{\mbox{\tinybox}}}

\newcommand{\weight}{\mathop{\rm w}\nolimits}

\newcommand{\seti}{{\mathcal I}}
\newcommand{\setj}{{\mathcal J}}
\newcommand{\setk}{{\mathcal K}}
\newcommand{\setl}{{\mathcal L}}
\newcommand{\sett}{{\mathcal T}}

\newcommand{\sets}{{\mathcal S}}
\newcommand{\seth}{{\mathcal H}}

\makeatletter
\renewcommand{\@biblabel}[1]{#1.}
\makeatother

\newtheorem{lemma}{Lemma}
\newtheorem{theorem}{Theorem}
\newtheorem{example}{Example}

\begin{document}
\begin{sloppypar}

\title{On the representation of bent functions 
by~bent~rectangles\footnote{
Appeared in
{\it Probabilistic Methods in Discrete Mathematics: 
Proceedings of the Fifth International Petrozavodsk Conference 
(Petrozavodsk, June 1--6, 2000)}. 
Utrecht, Boston: VSP, pp.~121--135, 2002.
}}

\author{Sergey Agievich\\
\small National Research Center
for Applied Problems of Mathematics and Informatics\\[-0.8ex]
\small Belarusian State University\\[-0.8ex]
\small Fr. Skorina av. 4, 220050 Minsk, Belarus\\[-0.8ex]
\small \texttt{agievich@bsu.by}
}

\date{}

\maketitle

\begin{abstract}
We propose a representation of boolean bent functions by bent
rectangles, that is, by special matrices with restrictions on rows and
columns. Using this representation, we
exhibit new classes of bent functions,
give an algorithm to construct bent functions,
improve a lower bound for the number of bent functions.
\end{abstract}

\section{Preliminaries}\label{Bent.1}

Let $V_n$ be an~$n$-dimensional vector space over the field
${\mathbb F}_2=\{0,1\}$
and $\func{F}_n$ be the set of all Boolean functions 
$V_n\to{\mathbb F}_2$.
We identify a function~$f\in\func{F}_n$ of~$x=(x_1,\ldots,x_n)$ 
with its algebraic normal form, that is, 
a polynomial of the ring~$\F_2[x_1,\ldots,x_n]$
reduced modulo the ideal~$(x_1^2-x_1,\ldots,x_n^2-x_n)$.
Denote by $\deg f$ the degree of such polynomial.

Write $\circup\gamma=(-1)^\gamma$ for $\gamma\in{\mathbb F}_2$.
The Walsh--Hadamard transform of $f$ is defined as
$$
\hatup{f}(\lambda)=
\sum_{x\in V_n}\circup f(x)(-1)^{\langle\lambda,x\rangle},\quad
\lambda\in V_n,
$$
where
$\langle\lambda,x\rangle=\lambda_1 x_1+\ldots+\lambda_n x_n$.
The symbol $+$ denotes both the addition in ${\mathbb F}_2$ and
the addition of integers or real numbers.
The way of the symbol usage depends on operands.

A function $f\in\func{F}_{2n}$ is called a {\it bent function} if
$|\hatup{f}(\lambda)|=2^n$ for any $\lambda\in V_{2n}$.
Let $\func{B}_{2n}$ be the set of all bent functions of $2n$ variables.

The Sylvester--Hadamard matrix $H_n$ of order $2^n$ is defined by the
recursive rule
$$
H_n=\left(
\begin{array}{rr}
1 &  1 \\
1 & -1 \\
\end{array}
\right)
\otimes H_{n-1},\quad
H_0=\left(1\right),
$$
where $\otimes$ denotes the Kronecker product.
Note that $H_n$ is symmetric and $H_n H_n=2^n E$,
where $E$ is the identity $2^n\times 2^n$ matrix.

Let $\func{A}_n=\{l\in\func{F}_n\colon \deg l\leq 1\}$ denote the set of
affine functions. A function $l\in\func{A}_n$ is of the form
$$
l(x)=\alpha_1 x_1+\ldots+\alpha_n x_n+\gamma=
\langle \alpha,x\rangle+\gamma,\quad
\alpha\in V_n,\quad
\gamma\in {\mathbb F}_2,
$$
and is a linear function if $\gamma=0$.

Let $v_1={\bf 0},v_2,\ldots,v_{2^n}$ be the
lexicographically ordered vectors of $V_n$.
The {\it sequence} and the {\it spectral sequence} of~$f$
are defined, respectively, as the row vectors
$$
\circup{\mathbf f}=(\circup f(v_1),\ldots,\circup f(v_{2^n}))
\quad\text{and}\quad
\hatup{\mathbf f}=(\hatup{f}(v_1),\ldots,\hatup{f}(v_{2^n})).
$$
It is known that the $i$th row of $H_n$ is a sequence of the
linear function $\langle\alpha, x\rangle$,
where $\alpha=v_i$
($\alpha$ is the $i$th vector of $V_n$)
and the Walsh--Hadamard transform can be written as
$$
\hatup{\mathbf f}=\circup{\mathbf f}H_n.
$$
Therefore,
\begin{itemize}
\item[(i)]
$\sum_{\lambda\in V_n}\hatup{f}^2 (\lambda)=2^{2n}$
(Parseval's identity);
\item[(ii)]
if
$l\in\func{A}_n$, $l(x)=\langle\alpha, x\rangle+\gamma$,
then
$\hatup{l}(\alpha)=\circup\gamma\cdot 2^n$
and
$\hatup{l}(\lambda)=0$, $\lambda\neq\alpha$.
\end{itemize}


\begin{lemma}\label{Lemma.BENT.SpectrumSum}
Let $f_1,f_2,f_3,f_4\in\func{F}_n$.
The vector
$
\frac{1}{2}
(
\hatup{\mathbf f}_1 + \hatup{\mathbf f}_2 +
\hatup{\mathbf f}_3 + \hatup{\mathbf f}_4
)$
is a spectral sequence of some function $g\in\func{F}_n$
if and only if
$$
f_1(x)+f_2(x)+f_3(x)+f_4 (x) =1
$$
for all $x\in V_n$.
Moreover,
$$
g(x)=f_1(x)f_2(x)+f_1(x)f_3(x)+f_2(x)f_3(x).
$$
\end{lemma}

\begin{proof}
If
$
\hatup{\mathbf g}=
\frac{1}{2}
(
\hatup{\mathbf f}_1+
\hatup{\mathbf f}_2+
\hatup{\mathbf f}_3+
\hatup{\mathbf f}_4
)
$, then
$
\circup{\mathbf g}=
\frac{1}{2}
(
\circup{\mathbf f}_1+
\circup{\mathbf f}_2+
\circup{\mathbf f}_3+
\circup{\mathbf f}_4
)
$
or
$$
\circup g(x)=
\frac{1}{2}
(
\circup f_1(x)+
\circup f_2(x)+
\circup f_3(x)+
\circup f_4(x)
),\quad
x\in V_n.
$$
But
$
\circup f_1(x)+
\circup f_2(x)+
\circup f_3(x)+
\circup f_4(x)=\pm 2$ if and only if
$$
f_1(x)+
f_2(x)+
f_3(x)+
f_4(x)=1.
$$

Let $g$ be the function defined in the statement of the lemma.
For all $x$, among $f_1(x)$, $f_2(x)$, $f_3(x)$
there are two identical values.
Without loss of generality we can assume that
$f_1(x)=f_2 (x)$.
Then
$g(x)=f_1(x)$, $f_4(x)=f_3(x)+1$,
and
\begin{align*}
\frac{1}{2}
(
\circup f_1(x)+
\circup f_2(x)+
\circup f_3(x)+
\circup f_4(x)
)&=
\frac{1}{2}
(
\circup f_1(x)+
\circup f_1(x)+
\circup f_3(x)-
\circup f_3(x)
)\\&=
\circup f_1(x)=
\circup g(x).
\end{align*}
Thus
$
\circup {\mathbf g}=
\frac{1}{2}
(
\circup{\mathbf f}_1+
\circup{\mathbf f}_2+
\circup{\mathbf f}_3+
\circup{\mathbf f}_4
)
$
and
$
\hatup{\mathbf g}=
\frac{1}{2}
(
\hatup{\mathbf f}_1+
\hatup{\mathbf f}_2+
\hatup{\mathbf f}_3+
\hatup{\mathbf f}_4
)$.
\end{proof}

Further we will use functions from the set
$$
\func{Q}_n=\{q\in\func{F}_n\colon \weight(\hatup{q})=4\},
$$
where~$\weight(\hatup{q})$ is the number of nonzero values of $\hatup{q}$.
The following lemma completely characterizes the elements of $\func{Q}_n$.

\begin{lemma}\label{Lemma.BENT.QSpectrum}
Let $q\in\func{Q}_n$ be such that
$\hatup{q}(\lambda)\neq 0$ only if
$$
\lambda\in\{
\lambda_{(1)},\lambda_{(2)},\lambda_{(3)},\lambda_{(4)}
\}\subseteq V_n,\quad
n\geq 2.
$$
Then
$\sum_i\lambda_{(i)}={\mathbf 0}$
and
$\hatup{q}(\lambda_{(i)})=\circup\gamma_i\cdot 2^{n-1}$
for some
$\gamma_i\in {\mathbb F}_2$ such that
$\sum_i\gamma_i=1$ or
$\prod_i\circup\gamma_i=-1$, $i=1,2,3,4$.
Moreover,
$$
q(x)=l_1(x)l_2(x)+l_1(x)l_3(x)+l_2(x)l_3(x),
$$
where
$l_i(x)=\langle\lambda_{(i)},x\rangle+\gamma_i$.
\end{lemma}

\begin{proof}
By induction on $n$ it is easy to show that the equation
$\sum_i z_i^2=2^{2n}$
has an unique solution $z_i=2^{n-1}$ in positive integers
(here and throughout the proof, $i$ varies from $1$ up to $4$).
According to the Parseval identity
$$
\sum_{i}\hatup{q}^2(\lambda_{(i)})=2^{2n}
$$
and $\hatup{q}(\lambda_{(i)})=\circup\gamma_i\cdot 2^{n-1}$
for some
$\gamma_i\in {\mathbb F}_2$.
But it yields
$
\hatup{\mathbf q}=\frac{1}{2}
(
\hatup{\mathbf l}_1+
\hatup{\mathbf l}_2+
\hatup{\mathbf l}_3+
\hatup{\mathbf l}_4
)
$,
where
$l_i(x)=\langle\lambda_{(i)},x\rangle+\gamma_i$.
The application of Lemma~\ref{Lemma.BENT.SpectrumSum} completes the proof.
\end{proof}

\section{Bent rectangles}\label{Bent.2}

A $2^m\times 2^k$ matrix with rows $a_{(i)}$
and transposed columns $a^{(j)}$
is called an {\it $(m,k)$ bent rectangle}
if each of the vectors $a_{(i)}$ and $2^{m-n}a^{(j)}$ is a spectral
sequence. Let $\boxup{\func{B}}_{m,k}$ be the set of all $(m,k)$ bent
rectangles.

\begin{theorem}\label{Theorem.BENT.Rect}
For all positive integers $m$, $k$ and $n$, $m+k=2n$, there exists
a bijection $\seth_{m,k}\colon \func{B}_{2n}\to\boxup{\func{B}}_{m,k}$.
\end{theorem}

\begin{proof}
Consider a function $f\in\func{B}_{2n}$ with the sequence
$$
\circup{\mathbf f}=(\circup{\mathbf f}_{(1)},\ldots,
\circup{\mathbf f}_{(M)}),\quad
\circup{\mathbf f}_{(i)}\in\{-1,1\}^K,\quad
M=2^m,\quad
K=2^k.
$$
Let $A$ be the $M\times K$ matrix with rows
\begin{equation}\label{Eq.BENT.Rect.1}
a_{(i)}=\circup{\mathbf f}_{(i)} H_k=\hatup{\mathbf f}_{(i)}.
\end{equation}
Assuming $H_m=(h_{ij})$, we obtain
\begin{multline*}
\circup{\mathbf f} H_{2n}
=\circup{\mathbf f}(H_m\otimes H_k)=\circup{\mathbf f}\left(
\begin{array}{rrr}
h_{11} H_k & \dotfill & h_{1M} H_k\\
\dotfill   & \ldots   & \dotfill\\
h_{M1} H_k & \dotfill & h_{MM} H_k\\
\end{array}
\right)=\\
=(
h_{11}a_{(1)}+\ldots+h_{M1} a_{(M)},
\ldots,
h_{1M} a_{(1)} +\ldots+h_{MM} a_{(M)}
)\in\{\pm 2^n\}^{2^{2n}}.
\end{multline*}
Hence
$a^{(j)}H_m\in\{\pm 2^n\}^{M}$,
where $a^{(j)}$ is the row vector that consists
of $j$th cooordinates of the vectors~$a_{(1)},\ldots,a_{(M)}$ or the
$j$th transposed column of $A$.
Therefore there exists a function $g\in\func{F}_m$ such that
$a^{(j)} H_m=2^n\circup{\mathbf g}$
or
\begin{equation}\label{Eq.BENT.Rect.2}
2^{m-n} a^{(j)}=2^{-n} a^{(j)} H_m H_m=
\circup{\mathbf g} H_m=\hatup{\mathbf g}.
\end{equation}
By~\eqref{Eq.BENT.Rect.1}, \eqref{Eq.BENT.Rect.2} the matrix
$A$ is an element of $\boxup{\func{B}}_{m,k}$.
Thus, we construct the map
$$
\seth_{m,k}\colon\quad
\func{B}_{2n}\to\boxup{\func{B}}_{m,k},\quad
f\mapsto
\left(
\begin{array}{c}
\circup{\mathbf f}_{(1)} H_k\\
\dotfill\\
\circup{\mathbf f}_{(M)} H_k\\
\end{array}
\right).
$$
It is easy to show that such map is bijective.
\end{proof}

Further we will denote
by~$\boxup f$ the image of~$f\in\func{B}_{2n}$ under the action
of the map $\seth_{m, k}$ constructed in the proof.
By construction, an $i$th row
of $\boxup f$ is a spectral sequence of the function
$$
f_{(i)}(x_{m+1},\ldots, x_{2n})=
f(\alpha_1,\ldots,\alpha_m,x_{m+1},\ldots,x_{2n})
$$
of $k$ variables, where
$\alpha=(\alpha_1,\ldots,\alpha_m)$
is the $i$th vector of $V_m$.


The obtained representation of bent functions by bent rectangles provides
a general way to receive many known results.
Consider, for example, the following construction due to Rothaus
(see~\cite[Class~II]{Rot76}).

\begin{example}
Let $f_1,f_2,f_3,f_4\in\func{B}_{2n}$ and 
\begin{align*}
f^*(y,z,x)&=f_1(x)f_2(x)+f_1(x)f_3(x)+f_2(x)f_3(x)\\
&+y(f_1(x)+f_2(x))+z(f_1(x)+f_3(x))+yz,\quad
x\in V_{2n},\quad
y,z\in{\mathbb F}_2.
\end{align*}
If~$f_1(x)+f_2(x)+f_3(x)+f_4(x)=0$ for all $x$,
then $f^*\in\func{B}_{2n+2}$. 

Indeed, extend~$\seth_{2,2n}$ to $\func{F}_{2n+2}$
and set~$\boxup f^*=\seth_{2,2n}(f^*)$.
By Lemma~\ref{Lemma.BENT.SpectrumSum},
$$
\boxup f^*=\frac{1}{2}
\left(
\begin{array}{c}
\hatup{\mathbf f}_1+
\hatup{\mathbf f}_2+
\hatup{\mathbf f}_3-
\hatup{\mathbf f}_4\\
\hatup{\mathbf f}_1-
\hatup{\mathbf f}_2+
\hatup{\mathbf f}_3+
\hatup{\mathbf f}_4\\
\hatup{\mathbf f}_1+
\hatup{\mathbf f}_2-
\hatup{\mathbf f}_3+
\hatup{\mathbf f}_4\\
\hatup{\mathbf f}_1-
\hatup{\mathbf f}_2-
\hatup{\mathbf f}_3-
\hatup{\mathbf f}_4\\
\end{array}
\right)
$$
and the transposed columns of $\boxup f^*$ have the form
$$
\frac{1}{2}
\left(
\hatup{f}_1(\lambda),
\hatup{f}_2(\lambda),
\hatup{f}_3(\lambda),
-\hatup{f}_4(\lambda)
\right)
H_2,
\quad
\lambda\in V_{2n}.
$$
But $|\hatup{f}_i(\lambda)|=2^n$ and,
multiplying such columns by $ 2^{1-n}$, we get spectral sequences.
Therefore $\boxup f^*\in\boxup{\func{B}}_{2,2n}$ and $f^*\in\func{B}_{2n+2}$.
\end{example}


Note that in~\cite{AdaTav90,SebZha94}
the $2^m\times 2^k$
matrices whose rows and columns are sequences of
bent functions were considered.
The multiplication by $2^{k/2}$ reduces such matrices
to the elements of $\boxup{\func{B}}_{m,k}$.

\section{Bent squares}\label{Bent.3}

We call~$\boxup f\in\boxup{\func{B}}_{n,n}$ a {\it bent square} and
write $\boxup{\func{B}}_n$ instead of $\boxup{\func{B}}_{n,n}$.

Suppose that the $i$th row (column) of $\boxup f$ is indexed by the $i$th
vector of $V_n$.
If $\boxup f\in\boxup{\func{B}}_n$ and $\seti,\setj$ are subsets of $V_n$, 
then $\boxup f[\seti\mid\setj]$ is the submatrix of $\boxup f$ determined by the
rows indexed by $\seti$ and columns indexed by $\setj$;
the submatrix $\boxup f(\seti\mid\setj)$ is obtained by removing such rows
and columns.

Let $\boxup{\func{B}}_n^{(\func{A})}$ 
be the set of all bent squares whose lines (rows
and columns) are spectral sequences of affine functions of $n$
variables.
A bent square $\boxup f\in\boxup{\func{B}}_n^{(\func{A})}$
have exactly one nonzero element $\pm 2^n$
in each line and the corresponding bent function $f$ has the form
\begin{equation}\label{Eq.BENT.AFunc}
f(x_1,\ldots,x_{2n})=
\varphi_1(x_1,\ldots,x_n)x_{n+1}+
\ldots
+\varphi_n(x_1,\ldots,x_n)x_{2n}+
\psi(x_1,\ldots,x_n),
\end{equation}
where~$\varphi_1,\ldots,\varphi_n,\psi\in\func{F}_n$, 
and the vector function 
$\varphi=(\varphi_1,\ldots,\varphi_n)$ determines a bijection on
$V_n$.
Indeed, for such $f$
$$
\boxup f[\{\alpha\}|\{\beta\}]=
\left\{
\begin{array}{rl}
\circup \psi(\alpha)\cdot 2^{n}, & \beta=\varphi(\alpha),\\
                              0 & \text{otherwise}.\\
\end{array}
\right.
$$

Functions of the form~\eqref{Eq.BENT.AFunc} constitute a
well-known Maiorana--McFarland class of bent functions~\cite{Dil72}.
Obviously,
$$
|\func{B}_{2n}^{(\func{A})}|=
|\boxup{\func{B}}_n^{(\func{A})}|=
2^{2^n}\cdot(2^n)!.
$$

Now
consider the set of bent squares $\boxup{\func{B}}_n^{(\func{AQ})}$ 
whose rows and,
consequently, columns are spectral sequences of functions from
$\func{A}_n\cup\func{Q}_n$.
Let~$\boxup f\in\boxup{\func{B}}_n^{(\func{AQ})}$ and let $\seti_f$
(respectively $\setj_f$) be the set of indices of rows
(respectively, columns) that are
spectral sequences of the elements of $\func{Q}_n$.
It is clear that
$|\seti_f|=|\setj_f|$.
For fixed $d=0,1,\ldots,2^n$, we define
$$
\boxup{\func{B}}_n^{(\func{AQ},d)}=
\{\boxup f\in\boxup{\func{B}}_n^{(\func{AQ})}
\colon
|\seti_f|=d
\}
$$
and let $\func{B}_{2n}^{(\func{AQ},d)}=
\seth_{n, n}^{-1}(\boxup{\func{B}}_n^{(\func{AQ},d)})$.
It is easy to see that $\boxup{\func{B}}_n^{(\func{AQ},d)}=\varnothing$
for $d=1,2,3,5$, $d\leq 2^n$.

We suggest a way to construct the function
$f^*\in\func{B}_{2n}^{(\func{AQ},4)}$ 
by using~$f\in\func{B}_{2n}^{(\func{A})}$, 
$n\geq 2$.
Suppose that $f$ is written in the form~\eqref{Eq.BENT.AFunc}
so that for some distinct $\alpha_{(i)}\in V_n$ it holds that
$$
\sum_i\alpha_{(i)}=\sum_i\beta_{(i)}={\mathbf 0},\quad
\beta_{(i)}=\varphi(\alpha_{(i)}),\quad
i=1,2,3,4.
$$
Next
\begin{itemize}
\item[(i)]
choose arbitrary $\gamma_{1j},\gamma_{2j},\gamma_{3j}\in
{\mathbb F}_2$,
define
$\gamma_{4j}=\gamma_{1j}+\gamma_{2j}+\gamma_{3j}+1$,
$j=1,2,3 $,
and set
$\gamma_{i4}=\gamma_{i1}+\gamma_{i2}+\gamma_{i3}+1$,
$i=1,2,3,4$;

\item[(ii)]
denote
$x_{(1)}=(x_1,\ldots,x_n)$,
$x_{(2)}=(x_{n+1},\ldots,x_{2n})$,
and for $v\in V_n$ put
$$
I_v(x_{(1)})=\prod_{i=1}^n (x_i+v_i+1)=
\left\{
\begin{array}{rl}
1, & x_{(1)}=v,\\
0 & \text{otherwise};
\end{array}
\right.
$$

\item[(iii)]
for $i=1,2,3,4$ define polynomials
$$
q_i(x_{(2)})=
f(\alpha_{(i)},x_{(2)})+
l_{i1}(x_{(2)})l_{i2}(x_{(2)})+
l_{i1}(x_{(2)})l_{i3}(x_{(2)})+
l_{i2}(x_{(2)})l_{i3}(x_{(2)}),
$$
where $l_{ij}(x_{(2)})=\langle\beta_{(j)},x_{(2)}\rangle+\gamma_{ij}$.
\end{itemize}
Now define $f^*$ by
\begin{equation}\label{Eq.BENT.AQ4Func}
f^*(x_1,\ldots,x_{2n})=
f^*(x_{(1)},x_{(2)})=f(x_{(1)},x_{(2)})+
\prod_{i=1}^4
I_{\alpha_{(i)}}(x_{(1)})
q_i(x_{(2)}).
\end{equation}
By the construction $f^*\in\func{B}_{2n}^{(\func{AQ},4)}$.
Indeed,
$$
\boxup f^*[\{\alpha_{(i)}\}|\{\beta_{(j)}\}]=\circup
\gamma_{ij}\cdot 2^{n-1}
$$
and
$$
\boxup f^*(\seti_{f^*}|\setj_{f^*})=\boxup f(\seti_{f^*}|\setj_{f^*}),\quad
\seti_{f^*}=\{\alpha_{(1)},\alpha_{(2)},\alpha_{(3)},\alpha_{(4)}\},\quad
\setj_{f^*}=\{\beta_{(1)},\beta_{(2)},\beta_{(3)},\beta_{(4)}\}.
$$

\begin{example}\label{Example.BENT.AQ4Func}
Consider the bent function
$$
f(x_1,x_2,x_3,x_4,x_5,x_6)=
x_1 x_2 x_3 +x_1 x_2+x_1 x_4 +x_2 x_6+x_3 x_5
$$
that can be represented as~\eqref{Eq.BENT.AFunc} with
$\varphi(x_1,x_2,x_3)=(x_1,x_3,x_2)$,
$\psi(x_1,x_2,x_3)=x_1 x_2 x_3 + x_1 x_2$.

Take
$$
\alpha_{(1)}=(0,0,0),\quad
\alpha_{(2)}=(0,0,1),\quad
\alpha_{(3)}=(1,0,0),\quad
\alpha_{(4)}=(1,0,1)
$$
and define the vectors
$$
\beta_{(1)}=(0,0,0),\quad
\beta_{(2)}=(0,1,0),\quad
\beta_{(3)}=(1,0,0),\quad
\beta_{(4)}=(1,1,0).
$$
Set
\begin{align*}
\gamma_{11}&=
\gamma_{12}=
\gamma_{13}=
\gamma_{21}=
\gamma_{22}=
\gamma_{31}=
\gamma_{33}=0,\\
\gamma_{23}&=\gamma_{32}=1.
\end{align*}
Now by~\eqref{Eq.BENT.AQ4Func}
$$
f^*(x_1,x_2,x_3,x_4,x_5,x_6)=
x_1 x_2 x_3 +
x_2 x_4 x_5 +
x_1 x_2 +
x_1 x_4 +
x_2 x_6 +
x_3 x_5 +
x_4 x_5
$$
and the bent squares $\boxup f$, $\boxup f^*$
have respectively the forms
$$
\left(
\begin{matrix}
\q 8&\q 0&\q 0&\q 0&\q 0&\q 0&\q 0&\q 0\\
\q 0&\q 0&\q 8&\q 0&\q 0&\q 0&\q 0&\q 0\\
\q 0&\q 8&\q 0&\q 0&\q 0&\q 0&\q 0&\q 0\\
\q 0&\q 0&\q 0&\q 8&\q 0&\q 0&\q 0&\q 0\\
\q 0&\q 0&\q 0&\q 0&\q 8&\q 0&\q 0&\q 0\\
\q 0&\q 0&\q 0&\q 0&\q 0&\q 0&\q 8&\q 0\\
\q 0&\q 0&\q 0&\q 0&\q 0&  -8&\q 0&\q 0\\
\q 0&\q 0&\q 0&\q 0&\q 0&\q 0&\q 0&\q 8\\
\end{matrix}
\right),\quad
\left(
\begin{matrix}
\q 4&\q 0&\q 4&\q 0&\q 4&\q 0&  -4&\q 0\\
\q 4&\q 0&\q 4&\q 0&  -4&\q 0&\q 4&\q 0\\
\q 0&\q 8&\q 0&\q 0&\q 0&\q 0&\q 0&\q 0\\
\q 0&\q 0&\q 0&\q 8&\q 0&\q 0&\q 0&\q 0\\
\q 4&\q 0&  -4&\q 0&\q 4&\q 0&\q 4&\q 0\\
  -4&\q 0&\q 4&\q 0&\q 4&\q 0&\q 4&\q 0\\
\q 0&\q 0&\q 0&\q 0&\q 0&  -8&\q 0&\q 0\\
\q 0&\q 0&\q 0&\q 0&\q 0&\q 0&\q 0&\q 8\\
\end{matrix}
\right).
$$
\end{example}

In the above construction we make the transition from
$f\in\func{B}_{2n}^{(\func{A})}$ to $f^*\in\func{B}_{2n}^{(\func{AQ},4)}$
based on the structure of the corresponding bent squares.
The following algorithm is intended for direct generation
of bent squares
$\boxup f\in\boxup{\func{B}}_n^{(\func{AQ},2d)}$,
$2\leq d\leq 2^{n-1}$,
with subsequent determination of corresponding bent functions.

The algorithm consists of the following steps:

\begin{enumerate}
\item
Choose sets
$\{\seti_1,\ldots,\seti_d\}$ and
$\{\setj_1,\ldots,\setj_d\}$
whose elements are disjoint $2$-element subsets of $V_n$ such that
$$
\sum_{\alpha\in\seti_1}\alpha=\ldots=\sum_{\alpha\in\seti_d}\alpha,\quad
\sum_{\beta\in\setj_1}\beta=\ldots=\sum_{\beta\in\setj_d}\beta.
$$
Let
$\seti=\seti_1\cup\ldots\cup\seti_d$,
$\setj=\setj_1\cup\ldots\cup\setj_d$.

\item
Choose a $0$-$1$ matrix $T=(t_{ij})$ of order $d$
with exactly $2$ ones in each line.
For $i,j=1,\ldots,d$, put
$$
\boxup f[\seti_i\mid\setj_j]=\left\{
\begin{array}{rl}
2^{n-1}
\left(
\begin{matrix}
\pm 1 & \pm 1\\
\pm 1 & \pm 1\\
\end{matrix}
\right),&
\mbox{if $t_{ij}=1$},\\
\left(
\begin{matrix}
0 & 0\\
0 & 0\\
\end{matrix}
\right),&
\mbox{if $t_{ij}=0$},\\
\end{array}
\right.
$$
where the signs are placed so that the product of nonzero elements in
each line of the matrix $2^{1-n}\boxup f[\seti\mid\setj]$ is equal
to~$-1$.

\item
Choose a matrix $\boxup f(\seti\mid\setj)$ so that in each its line there
is exactly one nonzero element~$\pm 2^n$.

\item
Put
$\boxup f[\{\alpha\}\mid\{\beta\}]=0$
if
$\alpha\in\seti$, $\beta\notin\setj$
or
$\alpha\notin\seti$, $\beta\in\setj$.
\end{enumerate}

By Lemma~\ref{Lemma.BENT.QSpectrum} the rows of $\boxup f$ indexed
by $\seti_f=\seti$ and the columns indexed by $\setj_f=\setj$
are spectral sequences of elements of $\func{Q}_n$.
The remaining lines are spectral sequences of affine functions,
$|\seti|=|\setj|=2d$,
and $\boxup f\in\boxup{\func{B}}_n^{(\func{AQ},2d)}$.

\begin{theorem}\label{Theorem.BENT.AQ2dCount}
The above algorithm allows to construct
$$
\left(
\left(\frac{(2^n-1)(2^{n-1})!}{(2^{n-1}-d)!}\right)^2\cdot
2^{4d}\cdot
\sum_{i=0}^d
\frac{(-1)^i}{i!}-
\delta_n(d)
\right)\cdot
(2^n-2d)!\cdot
2^{2^n-2d}
$$
distinct bent squares of the set $\boxup{\func{B}}_n^{(\func{AQ},2d)}$ 
for $d=2,\ldots,2^{n-1}$.
Here
$$
\delta_n(d)=
(d/2)!\cdot
2^{9d/2}\cdot
\left(
\frac{4}{3}s_n(d)r_n(d)-
\frac{4}{9}
r_n^2(d)
\right),
$$
if $d$ is even and 0 otherwise, and
$$
s_n(d)=(2^n-1)\cdot
\frac{(2^{n-1})!}{2^{d/2}(2^{n-1}-d)!(d/2)!},\quad
r_n(d)=(2^n-1)(2^{n-1}-1)\binom{2^{n-2}}{d/2}.
$$
\end{theorem}

Theorem~\ref{Theorem.BENT.AQ2dCount} gives a lower bound for
$|\boxup{\func{B}}_n^{(\func{AQ},2d)}|$ which is exact for~$d=2,3$.
Using the theorem, we obtain
$$
|\func{B}_8|=
|\boxup{\func{B}}_4|\geq
|\boxup{\func{B}}_4^{(\func{A})}|+
\sum_{d=2}^8
|\boxup{\func{B}}_4^{(\func{AQ},2d)}|\geq
1559994535674013286400\approx
2^{70.4}.
$$
This is the best lower bound for $|\func{B}_8|$ known to the author.

\section{Classification of $\func{B}_6$}

In~\cite{Rot76} Rothaus exhibits $4$ classes on~$\func{B}_6$ such
that elements of a same class can be obtained from each other
by linear invertible transformations of variables and addition of affine
terms.

We give an alternative classification of $\func{B}_6$: Two functions are
equivalent if the appropriate bent squares can be obtained one from 
the other by changing signs of elements and permuting rows and columns.
The representatives of each of the~$8$ found equivalence
classes are listed below as pairs
``bent square $\boxup f$~--- bent function $f(x_1,\ldots,x_6)$''.

\begin{description}
\item[Class 1]
($\func{B}_6^{(\func{A})}$,
contains $2^{15}\cdot 3^2\cdot 5\cdot 7$ elements):
\begin{gather*}
\left(
\begin{matrix}
\q 8&\q 0&\q 0&\q 0&\q 0&\q 0&\q 0&\q 0\\
\q 0&\q 8&\q 0&\q 0&\q 0&\q 0&\q 0&\q 0\\
\q 0&\q 0&\q 8&\q 0&\q 0&\q 0&\q 0&\q 0\\
\q 0&\q 0&\q 0&\q 8&\q 0&\q 0&\q 0&\q 0\\
\q 0&\q 0&\q 0&\q 0&\q 8&\q 0&\q 0&\q 0\\
\q 0&\q 0&\q 0&\q 0&\q 0&\q 8&\q 0&\q 0\\
\q 0&\q 0&\q 0&\q 0&\q 0&\q 0&\q 8&\q 0\\
\q 0&\q 0&\q 0&\q 0&\q 0&\q 0&\q 0&\q 8\\
\end{matrix}
\right),\\
x_1 x_4+
x_2 x_5+
x_3 x_6.
\end{gather*}

\item [Class 2]
($\func{B}_6^{(\func{AQ}, 4)}$,
contains $2^{18} \cdot 3\cdot 7^2$ elements):
\begin{gather*}
\left(
\begin{matrix}
 -4&\q 4&\q 4&\q 4&\q 0&\q 0&\q 0&\q 0\\
\q 4& -4&\q 4&\q 4&\q 0&\q 0&\q 0&\q 0\\
\q 4&\q 4& -4&\q 4&\q 0&\q 0&\q 0&\q 0\\
\q 4&\q 4&\q 4& -4&\q 0&\q 0&\q 0&\q 0\\
\q 0&\q 0&\q 0&\q 0&\q 8&\q 0&\q 0&\q 0\\
\q 0&\q 0&\q 0&\q 0&\q 0&\q 8&\q 0&\q 0\\
\q 0&\q 0&\q 0&\q 0&\q 0&\q 0&\q 8&\q 0\\
\q 0&\q 0&\q 0&\q 0&\q 0&\q 0&\q 0&\q 8\\
\end{matrix}
\right),\\
x_1 x_5 x_6+
x_1 x_4+
x_1 x_5+
x_1 x_6+
x_2 x_5+
x_3 x_6+
x_5 x_6+
x_5+
x_6.
\end{gather*}

\item[Class 3]
($\func{B}_6^{(\func{AQ},6)}$,
contains
$2^{21}\cdot 3\cdot 7^2$ elements):
\begin{gather*}
\left(
\begin{matrix}
  -4&\q 4&\q 4&\q 4&\q 0&\q 0&\q 0&\q 0\\
\q 4&  -4&\q 4&\q 4&\q 0&\q 0&\q 0&\q 0\\
\q 0&\q 0&  -4&\q 4&\q 4&\q 4&\q 0&\q 0\\
\q 0&\q 0&\q 4&  -4&\q 4&\q 4&\q 0&\q 0\\
\q 4&\q 4&\q 0&\q 0&  -4&\q 4&\q 0&\q 0\\
\q 4&\q 4&\q 0&\q 0&\q 4&  -4&\q 0&\q 0\\
\q 0&\q 0&\q 0&\q 0&\q 0&\q 0&\q 8&\q 0\\
\q 0&\q 0&\q 0&\q 0&\q 0&\q 0&\q 0&\q 8\\
\end{matrix}
\right),\\
x_1 x_2 x_6+
x_1 x_4 x_6+
x_1 x_5 x_6+
x_2 x_4 x_6+
x_1 x_5+
x_2 x_4+
x_2 x_5+
x_3 x_6+
x_5 x_6+
x_5+
x_6.
\end{gather*}

\item[Class 4]
($\func{B}_6^{(\func{AQ}, 7)}$,
contains
$2^{25}\cdot 3\cdot 7$ elements):
\begin{gather*}
\left(
\begin{matrix}
  -4&\q 4&\q 4&\q 4&\q 0&\q 0&\q 0&\q 0\\
\q 4&  -4&\q 0&\q 0&\q 4&\q 4&\q 0&\q 0\\
\q 4&\q 0&  -4&\q 0&\q 4&\q 0&\q 4&\q 0\\
\q 4&\q 0&\q 0&  -4&\q 0&\q 4&\q 4&\q 0\\
\q 0&\q 4&\q 4&\q 0&\q 0&\q 4&  -4&\q 0\\
\q 0&\q 4&\q 0&\q 4&  -4&\q 0&\q 4&\q 0\\
\q 0&\q 0&\q 4&\q 4&\q 4&  -4&\q 0&\q 0\\
\q 0&\q 0&\q 0&\q 0&\q 0&\q 0&\q 0&\q 8\\
\end{matrix}
\right),
\end{gather*}
\begin{multline*}
x_1 x_2 x_4+
x_1 x_3 x_4+
x_1 x_3 x_6+
x_1 x_4 x_5+
x_1 x_4 x_6+
x_1 x_5 x_6+\\+
x_2 x_3 x_4+
x_2 x_3 x_5+
x_2 x_3 x_6+
x_2 x_4 x_5+
x_2 x_5 x_6+
x_3 x_4 x_6+
x_3 x_5 x_6+\\+
x_1 x_5+
x_2 x_4+
x_2 x_5+
x_2 x_6+
x_3 x_4+
x_3 x_5+
x_3 x_6+
x_5 x_6+
x_5+
x_6.
\end{multline*}

\item[Class 5]
(a subset of $\func{B}_6^{(\func{AQ},8)}$,
contains
$2^{19} \cdot 7^2$ elements):
\begin{gather*}
\left(
\begin{matrix}
  -4&\q 4&\q 4&\q 4&\q 0&\q 0&\q 0&\q 0\\
\q 4&  -4&\q 4&\q 4&\q 0&\q 0&\q 0&\q 0\\
\q 4&\q 4&  -4&\q 4&\q 0&\q 0&\q 0&\q 0\\
\q 4&\q 4&\q 4&  -4&\q 0&\q 0&\q 0&\q 0\\
\q 0&\q 0&\q 0&\q 0&  -4&\q 4&\q 4&\q 4\\
\q 0&\q 0&\q 0&\q 0&\q 4&  -4&\q 4&\q 4\\
\q 0&\q 0&\q 0&\q 0&\q 4&\q 4&  -4&\q 4\\
\q 0&\q 0&\q 0&\q 0&\q 4&\q 4&\q 4&  -4\\
\end{matrix}
\right),\\
x_1 x_4+
x_2 x_5+
x_3 x_6+
x_5 x_6+
x_5+
x_6.
\end{gather*}

\item[Class 6]
(a subset of $\func{B}_6^{(\func{AQ},8)}$,
contains
$2^{20}\cdot 3^2\cdot 7^2$ elements):
\begin{gather*}
\left(
\begin{matrix}
  -4&\q 4&\q 4&\q 4&\q 0&\q 0&\q 0&\q 0\\
\q 4&  -4&\q 4&\q 4&\q 0&\q 0&\q 0&\q 0\\
\q 0&\q 0&  -4&\q 4&\q 4&\q 4&\q 0&\q 0\\
\q 0&\q 0&\q 4&  -4&\q 4&\q 4&\q 0&\q 0\\
\q 0&\q 0&\q 0&\q 0&  -4&\q 4&\q 4&\q 4\\
\q 0&\q 0&\q 0&\q 0&\q 4&  -4&\q 4&\q 4\\
\q 4&\q 4&\q 0&\q 0&\q 0&\q 0&  -4&\q 4\\
\q 4&\q 4&\q 0&\q 0&\q 0&\q 0&\q 4&  -4\\
\end{matrix}
\right),\\
x_2 x_4 x_6+
x_1 x_4+
x_2 x_4+
x_2 x_5+
x_3 x_6+
x_5 x_6+
x_5+
x_6.
\end{gather*}

\item[Class 7]
(a subset of $\func{B}_6^{(\func{AQ},8)}$,
contains
$2^{23}\cdot 3\cdot 7^2$ elements):
\begin{gather*}
\left(
\begin{matrix}
  -4&\q 4&\q 4&\q 4&\q 0&\q 0&\q 0&\q 0\\
\q 4&  -4&\q 0&\q 0&\q 4&\q 4&\q 0&\q 0\\
\q 4&\q 0&  -4&\q 0&\q 4&\q 0&\q 4&\q 0\\
\q 4&\q 0&\q 0&  -4&\q 0&\q 4&\q 4&\q 0\\
\q 0&\q 4&\q 4&\q 0&  -4&\q 0&\q 0&\q 4\\
\q 0&\q 4&\q 0&\q 4&\q 0&  -4&\q 0&\q 4\\
\q 0&\q 0&\q 4&\q 4&\q 0&\q 0&  -4&\q 4\\
\q 0&\q 0&\q 0&\q 0&\q 4&\q 4&\q 4& -4\\
\end{matrix}
\right),
\end{gather*}
\begin{multline*}
x_1 x_2 x_4+
x_1 x_2 x_5+
x_1 x_2 x_6+
x_1 x_3 x_4+
x_1 x_3 x_5+
x_1 x_3 x_6+
x_1 x_4 x_5+\\+
x_1 x_4 x_6+
x_2 x_3 x_4+
x_2 x_3 x_5+
x_2 x_3 x_6+
x_2 x_4 x_5+
x_2 x_5 x_6+
x_3 x_4 x_6+
x_3 x_5 x_6+\\+
x_2 x_4+
x_2 x_6+
x_2 x_5+
x_3 x_4+
x_3 x_5+
x_3 x_6+
x_5 x_6+
x_5+
x_6.
\end{multline*}

\item[Class 8]
(
contains
$2^{23}\cdot 3^2\cdot 5\cdot 7$ elements):
\begin{gather*}
\left(
\begin{matrix}
  -6&\q 2&\q 2&\q 2&\q 2&\q 2&\q 2&\q 2\\
\q 2&  -6&\q 2&\q 2&\q 2&\q 2&\q 2&\q 2\\
\q 2&\q 2&  -6&\q 2&\q 2&\q 2&\q 2&\q 2\\
\q 2&\q 2&\q 2&  -6&\q 2&\q 2&\q 2&\q 2\\
\q 2&\q 2&\q 2&\q 2&  -6&\q 2&\q 2&\q 2\\
\q 2&\q 2&\q 2&\q 2&\q 2&  -6&\q 2&\q 2\\
\q 2&\q 2&\q 2&\q 2&\q 2&\q 2&  -6&\q 2\\
\q 2&\q 2&\q 2&\q 2&\q 2&\q 2&\q 2&  -6\\
\end{matrix}
\right),\\
x_4 x_5 x_6+
x_1 x_4+
x_2 x_5+
x_3 x_6+
x_4 x_5+
x_4 x_6+
x_5 x_6+
x_4+
x_5+
x_6.
\end{gather*}
\end{description}

\section{Proof of Theorem~\ref{Theorem.BENT.AQ2dCount}}

Consider a $0$-$1$ matrix $A$ without zero lines.
Denote by $\omega(A)$ a set of matrices obtained
by all arrangements of signs in $A$
such that the product of nonzero elements in each line
is equal to~$-1$.

Let $\sett$ denote the set of all $0$-$1$ matrices containing exactly
2 ones in each line and let $\sett^*$ be the set of all such
matrices without proper submatrices with the same property.
For $\sets\subseteq\sett$, define the set
$$
\Omega(\sets)=
\bigcup_{T\in\sets}
\omega\left(
T\otimes\left(
\begin{matrix}
1 & 1\\
1 & 1\\
\end{matrix}
\right)
\right)
$$
and its exponential generating function
$$
\Psi_{\Omega(\sets)}(x)=
\sum_{d\geq 2}|\Omega(\sets_d)|\cdot \frac{x^d}{(d!)^2}.
$$
Here $\sets_d$ is a set of all $d\times d$ matrices contained in $\sets$.

\begin{lemma}\label{Lemma.BENT.Coloring}
Let $A$ be a $0$-$1$ matrix of order $d$ with $m$ ones.
If there exists a matrix
$B\in\sett^*_d$ such that $A\geq B$,
then 
$$
|\omega(A)|=2^{m-2d+1}.
$$
\end{lemma}

\begin {proof}
Without loss of generality we can assume that
$$
B=\left(
\begin{matrix}
1 & 1 & 0 & 0 &\ldots & 0 & 0\\
0 & 1 & 1 & 0 &\ldots & 0 & 0\\
0 & 0 & 1 & 1 &\ldots & 0 & 0\\
\dots & \dots & \dots & \dots & \ldots & \dots & \dots\\
0 & 0 & 0 & 0 &\ldots & 1 & 1\\
1 & 0 & 0 & 0 &\ldots & 0 & 1\\
\end{matrix}
\right).
$$
Indeed, by a permutation of rows and columns
we can reduce any matrix of~$\sett^*_d$
to this form although such permutations
do not change the number $|\omega(A)|$.

Form a matrix $\bar C$ by assigning appropriate signs to the elements
of the $0$-$1$ matrix $C=A-B$. Let $\bar r_i$, $\bar s_j$ be the products,
respectively, of nonzero elements of the $i$th row and the $j$th column of
$\bar C$ (we suppose that such products are equal to $1$ for zero lines). The
matrix $C$ contains $m-2d$ ones and signs can be assigned in $2^{m-2d}$
ways.

If the matrix $\bar B=(\bar b_{ij})$ is obtained by an arrangement
of signs of ones in $B$ such that $\bar A=\bar B+\bar C\in\omega(A)$, then
\begin{equation}\label{Eq.BENT.System}
\bar b_{11}\bar b_{1d}=-\bar s_1,\
\bar b_{11}\bar b_{12}=-\bar r_1,\
\bar b_{12}\bar b_{22}=-\bar s_2,\
\bar b_{22}\bar b_{23}=-\bar r_2,\
\ldots,\
\bar b_{d1}\bar b_{dd}=-\bar r_d.
\end{equation}
Any nonzero element of~$\bar C$ enters once in each of the products
$\prod_{i=1}^d\bar r_i$ and
$\prod_{j=1}^d\bar s_j$,
therefore
$\prod\bar r_i\prod \bar s_j=1$
and there exist two distinct matrices $\bar B$ that meet~\eqref{Eq.BENT.System}.
Hence there are exactly
$2\cdot 2^{m-2d}$ distinct matrices
$\bar A\in\omega(A)$, as claimed.
\end {proof}

Let us determine the numbers $t(d)=|\Omega(\sett_d)|$, $d\geq 2$.
Observe that any matrix of $\Omega(\sett)$ by a permutation of rows
and columns can be represented as the direct sum of matrices contained in
$\Omega(\sett^*)$.
Using the exponential generating
functions theory~\cite[Chapter~3]{GouJac83},
we obtain
$$
\Psi_{\Omega(\sett)}(x)=
\exp\left(
\Psi_{\Omega(\sett^*)}(x)
\right)-1
$$
and
$$
t(d)=
\left[\frac{x^d}{(d!)^2}\right]
\Psi_{\Omega(\sett)}(x)=
\left[\frac{x^d}{(d!)^2}\right]
\exp\left(\Psi_{\Omega(\sett^*)}(x)\right).
$$

It is known~\cite{GouJac83}, that 
$$
|\sett_d^*|=\frac{(d!)^2}{2d},\quad
d\geq 2.
$$
Besides, if
$$
A=T\otimes
\left(
\begin{matrix}
1 & 1\\
1 & 1\\
\end{matrix}
\right)
$$
and the matrix $T\in\sett^*$ is represented as the sum 
of the nonzero $0$-$1$
matrices $T_1$ and $T_2$, then the matrix
$$
B=
T_1\otimes
\left(
\begin{matrix}
1 & 0\\
0 & 1\\
\end{matrix}
\right)+
T_2\otimes
\left(
\begin{matrix}
0 & 1\\
1 & 0\\
\end{matrix}
\right)\in\sett_{2d}^*.
$$
Moreover $A\geq B$ and 
$|\omega(A)|=2^{4d+1}$
by Lemma~\ref{Lemma.BENT.Coloring}.
Thus,
$$
\Psi_{\Omega(\sett^*)}(x)=
\sum_{d\geq 2}
\frac{(d!)^2}{2d}\cdot
2^{4d+1}\cdot
\frac{x^d}{(d!)^2}=
\log(1-16x)^{-1}-16x
$$
and
$$
t(d)=
\left[\frac{x^d}{(d!)^2}\right]
\exp(-16x)(1-16x)^{-1}=
2^{4d}\cdot(d!)^2\cdot\sum_{i=0}^d\frac{(-1)^i}{i!}.
$$

Return to the analyzed algorithm. The output
bent squares $\boxup f\in\boxup B_n^{(AQ,2d)}$ are uniquely
determined by
\begin{enumerate}
\item[(i)]
the choice of each of the sets $\{\seti_1,\ldots,\seti_d\}$ and
$\{\setj_1,\ldots,\setj_d\}$ in $(2^n-1)\binom{2^{n-1}}{d}$ ways;

\item[(ii)]
the choice of the matrix
$\boxup f[\seti\mid\setj]$ in $t(d)$ ways;

\item[(iii)]
the choice of the matrix $\boxup f(\seti\mid\setj)$
in $2^{2^n-2d}\cdot(2^n-2d)!$ ways.
\end{enumerate}

For distinct input data the output bent squares are distinct
except for, possibly, the case:
$d$ is even and the matrix
$\boxup f[\seti\mid\setj]$ by a permutation of rows and
columns can be represented as the direct sum of matrices from the set
$\omega(P)$, 
where 
$$
P=\left(\begin{matrix}
1&1&1&1\\
1&1&1&1\\
1&1&1&1\\
1&1&1&1
\end{matrix}
\right).
$$
Let $S_n(d)$ be the number of possible sets of input
data for this case and $R_n(d)$ be the number of distinct 
bent squares obtained from such sets. 
Then the algorithm allows us to construct
\begin{equation}\label{Eq.BENT.Proof1}
\left(
(2^n-1)\binom{2^{n-1}}{d}
\right)^2\cdot
t(d)\cdot
2^{2^n-2d}\cdot
(2^n-2d)!-S_n(d)+R_n(d)
\end{equation}
distinct bent squares. 
It remains to obtain the numbers~$S_n(d)$, $R_n(d)$.

In the interesting case
for some choice of the sets
$\{\{\seti_1,\seti_2\},\ldots,\{\seti_{d-1},\seti_{d}\}\}$,
$\{\{\setj_1,\setj_2\},\ldots$ $\ldots,\{\setj_{d-1},\setj_{d}\}\}$
and a permutation~$\pi$ on $\{1,\ldots,d/2\}$, 
it holds that
$$
\boxup f[\setk_i\mid\setl_{\pi(i)}]=
\bar P_i,\quad
\bar P_i\in\omega(P),\quad
i=1,\ldots,d/2,
$$
where
$\setk_i=\seti_{2i-1}\cup\seti_{2i}$,
$\setl_i=\setj_{2i-1}\cup\setj_{2i}$.

Each of the sets
$\left\{
\{\seti_1,\seti_2\},\ldots,\{\seti_{d-1},\seti_{d}\}
\right\}$
and
$\left\{
\{\setj_1,\setj_2\},\ldots,\{\setj_{d-1},\setj_{d}\}
\right\}$
can be chosen in~$s_n(d)$ ways. Further, there are
$(d/2)!$ ways to choose the permutation
$\pi$ and
$2^9$ ways to choose each of the matrices $\bar P_i$.
Thus
\begin{equation}\label{Eq.BENT.Proof2}
S_n(d)=s_n^2(d)\cdot
(d/2)!\cdot
2^{9d/2}\cdot
2^{2^n-2d}\cdot
(2^n-2d)!.
\end{equation}

Let $L(\alpha,\beta)$ denote
the $2$-dimension linear subspace of $V_n$,
generated by vectors
$\alpha$ and $\beta$, and let
$\gamma_1,\ldots,\gamma_{d/2}$ be the
representatives of distinct classes
$V_n/L(\alpha,\beta)$.
The pair of the sets
$\{\alpha,\beta\}$,
$\{\gamma_1,\ldots,\gamma_{d/2}\}$
can be chosen in $r_n(d)$ ways.

If
$$
\seti_{2i-1}=\gamma_i+\{{\bf 0},\alpha\},\quad
\seti_{2i}=\gamma_i+\{\beta,\alpha+\beta\}
$$
or
$$
\seti_{2i-1}=\gamma_i+\{{\bf 0},\beta\},\quad
\seti_{2i}=\gamma_i+\{\alpha,\alpha+\beta\}
$$
or
$$
\seti_{2i-1}=\gamma_i+\{{\bf 0},\alpha+\beta\},\quad
\seti_{2i}=\gamma_i+\{\alpha,\beta\}
$$
for $i=1,\ldots,d/2$,
then
$$
\setk_i=\gamma_i+L(\alpha,\beta).
$$
Thus only ${r_n(d)}/{3}$
distinct sets $\{\setk_1,\ldots,\setk_{d/2}\}$
correspond to $r_n(d)$ distinct sets
$\{\{\seti_1,\seti_2\},\ldots$ $\ldots,\{\seti_{d-1},\seti_{d}\}\}$.
Similarly,
${r_n(d)}/{3}$ distinct sets $\{\setl_1,\ldots,\setl_{d/2}\}$
correspond to $r_n(d)$ distinct sets
$\{\{\setj_1,\setj_2\},\ldots,\{\setj_{d-1},\setj_{d}\}\}$.
Therefore,
\begin{multline}\label{Eq.BENT.Proof3}
R_n(d)=
\left(
(s_n(d)-r_n(d))^2+
2(s_n(d)-r_n(d))\frac{r_n(d)}{3}+
\left(\frac{r_n(d)}{3}\right)^2
\right)\times\\
\times(d/2)!\cdot
2^{9d/2}\cdot
2^{2^n-2d}\cdot
(2^n-2d)!.
\end{multline}

Substituting~\eqref{Eq.BENT.Proof2} and~\eqref{Eq.BENT.Proof3}
into~\eqref{Eq.BENT.Proof1}, we obtain the desired result.

\end{sloppypar}
\end{document}